\documentclass{aims}
\usepackage{comment}
\usepackage{amsmath,amssymb}
\usepackage{paralist}
\usepackage{graphics} 
\usepackage{epsfig} 
\usepackage{graphicx}  \usepackage{epstopdf}
\usepackage[colorlinks=true]{hyperref}
\hypersetup{urlcolor=blue, citecolor=red}

\def\currentvolume{X}
\def\currentissue{X}
\def\currentyear{200X}
\def\currentmonth{XX}
\def\ppages{X--XX}
\def\DOI{10.3934/xx.xx.xx.xx}

\theoremstyle{definition}

\title[] 
{Traveling wave and aggregation in a flux-limited Keller-Segel model}

\author[Vincent Calvez, Beno\^it Perthame, and Shugo Yasuda]{}

\subjclass{Primary: 92C17 35C07; Secondary: .}
\keywords{Flux-limited Keller-Segel model, Chemotaxis, Traveling wave, Aggregation, Population growth.}

\email{vincent.calvez@math.cnrs.fr}
\email{benoit.perthame@upmc.fr}
\email{yasuda@sim.u-hyogo.ac.jp}

\begin{document}
\maketitle

\centerline{\scshape Vincent Calvez}
\medskip
{\footnotesize
	\centerline{Institut Camille Jordan, UMR 5208 CNRS/Universit\'e Claude Bernard Lyon 1, and Project-team Inria NUMED, Lyon}
} 

\medskip

\centerline{\scshape Beno\^it Perthame}
\medskip
{\footnotesize
	\centerline{Sorbonne Universit\'es, UPMC Univ Paris 06, Laboratoire Jacques-Louis Lions UMR CNRS 7598,}
	\centerline{Universit\'e Paris Diderot, Inria de Paris, F75005 Paris, France}
} 

\medskip

\centerline{\scshape Shugo Yasuda}
\medskip
{\footnotesize
	\centerline{Graduate School of Simulation Studies, University of Hyogo, Kobe 650-0047, Japan}
}

\bigskip

\centerline{(Communicated by the associate editor name)}

\begin{abstract}
Flux-limited Keller-Segel (FLKS) model has been recently derived from kinetic transport models for bacterial chemotaxis and shown to represent better the collective movement observed experimentally.
Recently, associated to the kinetic model, a new instability formalism has been discovered related to stiff chemotactic response.
This motivates our study of traveling wave and aggregation in population dynamics of chemotactic cells based on the FLKS model with a population growth term.
	
Our study includes both numerical and theoretical contributions.
In the numerical part, we uncover a variety of solution types in the one-dimensional FLKS model additionally to standard Fisher/KPP type traveling wave.
The remarkable result is a counter-intuitive backward traveling wave, where the population density initially saturated in a stable state transits toward an unstable state in the local population dynamics. 
Unexpectedly, we also find that the backward traveling wave solution transits to a localized spiky solution as increasing the stiffness of chemotactic response.

In the theoretical part, we obtain a novel analytic formula for the minimum traveling speed which includes the counter-balancing effect of chemotactic drift vs. reproduction/diffusion in the propagating front.
The front propagation speeds of numerical results only slightly deviate from the minimum traveling speeds, except for the localized spiky solutions, even for the backward traveling waves.
We also discover an analytic solution of unimodal traveling wave in the large-stiffness limit, which is certainly unstable but exists in a certain range of parameters. 
\end{abstract}

\section{Introduction}
Aggregations and traveling waves are ubiquitous in collective dynamics of chemotactic cells. 
It is well known that chemotactic bacteria as \emph{E. Coli} extend their habitat as creating patterns with localized aggregations of population \cite{art:66A, art:91BB, art:95BB, art:11Letal}. 
Aggregation stems from the chemotaxis of motile cells, where cells are attracted to migrate toward a higher-concentration region of chemical cues produced by themselves.
A challenge is to understand the interaction between aggregation by chemotaxis and invasion by reproduction/diffusion as in the usual Fisher/KPP model.

A reaction-diffusion-advection equation, to describe the population dynamics of chemotactic cells was proposed by Keller and Segel \cite{art:71KS, art:71KS2}, where the reaction term describes the local population dynamics such as cell proliferation, the diffusion term describes the random motions of cells, and the advection term describes the chemotaxis flux of cell population.
Modification of the classical Keller-Segel (KS) model has also been developed to obtain a model which can describe more realistic behaviors of chemotactic cells \cite{art:08TMPA,art:09HP}.

Flux-limited Keller-Segel (FLKS) model, which is a very active research subject nowadays \cite{art:12CKWW, art:15BBTW, art:17BW}, can avoid nonphysical blow-up phenomena due to unbounded chemotaxis flux inherent in the classical KS model, and thus describes more realistically the collective dynamics of chemotactic cells. 
The boundedness of chemotaxis flux in the FLKS model is related to a biological function in chemotactic sensing of cells, i.e., the stiff and bounded signal response \cite{art:83BSB}. 
This can be seen in derivation of the FLKS model from a related kinetic chemotaxis model with stiff chemotactic response \cite{art:05DS, art:11SCBPBS, art:13JV, art:16C, art:17PVW}.

Numerical studies on traveling waves and aggregations in the KS model with a population growth term have been carried out for various biological systems, and thus the pattern formation mechanism of chemotactic cells and the mathematical properties of the spatio-temporal dynamics have been investigated \cite{art:08TMPA, art:91MMWM, art:96MT, art:11PH, art:14EIM}. 
The traveling pulses in the FLKS model have also been investigated numerically, and thus the importance of stiffness and modulation amplitude in the bounded chemotactic flux is clarified to reproduce the collective migrations of chemotactic cells \cite{art:10SCBBSP, art:16EGBAV}.

In this paper, we study the traveling wave and aggregation in the one-dimensional FLKS model with a population growth term both theoretically and numerically. 
In the numerical part, we put a focus on the effects of stiffness and modulation amplitude in chemotaxis flux and unveil a variety of solution types in the FLKS model according to the stiffness and modulation parameters. 
In the theoretical part, we analytically calculate the traveling speed in front propagation and the unimodal traveling wave solution in the stiff-flux limit.

\section{Basic equation and preliminary analysis}
\subsection{Basic equation}
We consider a one-dimensional FLKS equation,
\begin{equation}\label{eq_KS}
	\partial_t\rho(t,x)+\partial_x(U_\delta[\partial_x \log S(t,x)]\rho)=\partial_{xx}\rho + P[\rho]\rho,
\end{equation} 
where $\rho(t,x)$ is the population density of cells at position $x\in \mathbb{R}$ and time $t \ge 0$ and $S(t,x)$ is the concentration of chemoattractant.

The second term of the left-hand side of Eq.~(\ref{eq_KS}) represents a flux due to chemotaxis of cells, where $U_\delta(X)$ is a bounded increasing function written as
\begin{equation}\label{eq_U}  
	U_\delta(X)=U\left(\frac{X}{\delta}\right),\quad U'(X)>0,\quad
	U(X)\rightarrow \pm\chi\quad(X\rightarrow \pm\infty).
	\end{equation}
Here $\chi\,(>0)$ and $\delta^{-1}\,(>0)$ represent the modulation and stiffness in the chemotactic response of cells.

The second term of the right-hand side of Eq.~(\ref{eq_KS}) describes a proliferation of cells, where $P[\rho]$ represents a proliferation rate which we choose, for simplicity, as
\begin{equation}\label{eq_P}
	P(\rho)=\left\{
		\begin{array}{c c}
			p&(0\le\rho\le\rho_c),\\
			\frac{1}{\rho}-1&(\rho_c<\rho),
		\end{array}
		\right.
	\end{equation} 
with
\begin{equation}\label{eq_rhoc}
	\rho_c=\frac{1}{1+p}.
\end{equation}
The proliferation rate $P(\rho)$ is positive and constant, $p\,(>0)$, in a lower-density regime ($0<\rho<\rho_c$), but it monotonically decreases in a higher-density regime ($\rho>\rho_c$) and becomes negative for $\rho>1$ so that the population saturates to $\rho=1$ in the higher-density regime.
Here, the constant $p$ represents a relative amplitude of the proliferation rate in the lower-density regime to the rate of change toward the saturated state $\rho=1$ in the higher-density regime.

The concentration of chemoattractant $S(t,x)$ is produced by the chemotactic cells themselves and described by
\begin{equation}\label{eq_S} 
	-d\partial_{xx}S+S=\rho,
\end{equation}
where $d$ is the diffusion coefficient, $d>0$.

We remark that the boundedness of chemotaxis drift velocity in Eq.~(\ref{eq_U}) stems from a stiff and bounded signal response of chemotactic cells to the external chemoattractant concentration $S(t,x)$.
In the chemotaxis response we also consider the logarithmic sensing, where cells can sense a relative variation of chemoattractant concentration to the local one along their moving pathway \cite{art:09KJTW}.
These microscopic backgrounds involved in the chemotaxis drift velocity formalism can be described in a kinetic transport model at the individual level, and the FLKS equation Eq.~(\ref{eq_KS}) can be derived by the asymptotic analysis of the kinetic chemotaxis model with stiffness \cite{art:17PVW}.

It is easily seen that the above equations admit two constant stationary states at $\rho=S=0$ and $\rho=S=1$.
The steady state $\rho=S=0$ is unconditionally unstable while the steady state $\rho=S=1$ is conditionally linearly stable.
The linear stability condition of Eq.~(\ref{eq_KS}) is written as \cite{art:08NPR},
\begin{equation}\label{eq_stability}
	U_\delta'[0]\le (1+\sqrt{d})^2.
\end{equation}
We remark that the linear stability condition Eq.~(\ref{eq_stability}) is also obtained by a diffusion scaling of the linear instability condition obtained in the related kinetic chemotaxis model \cite{art:17PY}.
In Ref. \cite{art:17PY}, it is also confirmed that the bounded periodic patterns appear from the initial uniform state with $\rho=1$ when the linear stability condition Eq.~(\ref{eq_stability}) is violated, as phenomena similar to Turing instability because localized patterns occur.

In this paper, we are concerned with the solution which connects the unstable state at $\rho=S=0$ and the conditionally stable state at $\rho=S=1$ for the FLKS system.
In the following, we consider the boundary condition
\begin{equation}\label{eq_bc}
	\begin{array}{cc}
		\rho=S=1, &\mathrm{at}\quad x=-\infty,\\
		\rho=S=0, &\mathrm{at}\quad x=\infty.
	\end{array}
\end{equation}

\subsection{Traveling speed}\label{sec_analy_speed}

We introduce a coordinate $\xi$ relative to the moving frame with a constant velocity $c$, i.e., $\xi=x-ct$ and consider the Cauchy problem for the population density $\rho(t,x)=\tilde \rho(\xi)$ and concentration of chemoattractant $S(t,x)=\tilde S(\xi)$.
Then, we can rewrite Eqs. (\ref{eq_KS}) and (\ref{eq_S}) as
\begin{equation}\label{eq_trho}
	-c\tilde \rho'+(U_\delta[(\log \tilde S)']\tilde \rho)'=\tilde \rho''+P(\tilde \rho)\tilde \rho,
\end{equation} 
and
\begin{equation}\label{eq_cauchy_S}
	-d\tilde S''+\tilde S=\tilde \rho.
\end{equation}
Here, the prime $\prime$ represents the derivative with respect to $\xi$.
The boundary conditions are written as
\begin{equation}\label{eq_bc_trho}
	\begin{array}{cc}
		\tilde\rho(\xi)=\tilde S(\xi)=1&\xi\rightarrow -\infty,\\
		\tilde\rho(\xi)=\tilde S(\xi)=0&\xi\rightarrow \infty.
	\end{array}
\end{equation} 
Furthermore, we can write Eq.~(\ref{eq_cauchy_S}) in the convolution
\begin{equation}\label{eq_tS} 
	\tilde S(\xi)=\frac{1}{2\sqrt{d}}
	\int_{-\infty}^\infty e^{-\frac{|\xi-\zeta|}{\sqrt{d}}}
	\tilde\rho(\zeta)d\zeta.
\end{equation}

Following \cite{art:16C,book:15P,art:11NPT,art:37KPP}, one can establish formally a relation between the decay rate of population density at $\xi\gg 1$ and the traveling speed.
To do so, we consider the exponential decay of $\tilde \rho(\xi)$ at far-right region, i.e.,
\begin{equation}\label{eq_app_trho}
	\tilde \rho(\xi)\propto e^{-\lambda\xi}, \quad \xi \gg 1.
\end{equation}
Then, from Eq.~(\ref{eq_tS}), we can write
\begin{equation}\label{eq_app_tS}
	\tilde S(\xi)\propto e^{-\min (\lambda,\frac{1}{\sqrt{d}})\xi},\quad \xi \gg 1.
\end{equation}

By substituting Eqs.~(\ref{eq_app_trho}) and (\ref{eq_app_tS}) into Eq.~(\ref{eq_trho}), we obtain a formula for propagation speed $c$ as a function of decay rate $\lambda$, i.e.,
\begin{equation}\label{eq_cl} 
	c(\lambda)=\lambda+\frac{p}{\lambda}
	-U_\delta\left[\min \left(\lambda,\frac{1}{\sqrt{d}}\right)\right].
\end{equation}
We note that from Eq.~(\ref{eq_app_tS}), the derivative of $\log\tilde S$ is constant, i.e., $(\log \tilde S)'=-\min (\lambda,\frac{1}{\sqrt{d}})$, at $\xi\gg 1$.

Furthermore, for the following flux function,
\begin{equation}\label{eq_arcU} 
	U_\delta(X)=\frac{2\chi}{\pi}\arctan\left(\frac{X}{\delta}\right),
\end{equation}
we can calculate the minimum speed $c$ of Eq.~(\ref{eq_cl}) analytically as
\begin{equation}\label{eq_minc}
\min_\lambda c(\lambda)=\left\{
	\begin{array}{ccc}
		2\sqrt{p}-\frac{2\chi}{\pi}\arctan\left(\frac{1}{\sqrt{d}\delta}\right),&
		\mathrm{if}&dp>1,\\
		\frac{1}{\sqrt{d}}+p\sqrt{d}-\frac{2\chi}{\pi}\arctan\left(\frac{1}{\sqrt{d}\delta}\right),
		&\mathrm{if}&1-\frac{2\chi}{\pi(\delta+\frac{1}{d\delta})}<dp<1,\\
		\Lambda+\frac{p}{\Lambda}-\frac{2\chi}{\pi}\arctan\left(\frac{\Lambda}{\delta}\right),
		&\mathrm{if}&dp\le 1-\frac{2\chi}{\pi(\delta+\frac{1}{d\delta})}
	\end{array}
	\right.
\end{equation}
where $\Lambda$ is defined as
\begin{equation}\label{eq_Lambda}
	\Lambda=\left\{
		\frac{p-\delta^2+\frac{2\chi}{\pi}\delta+\sqrt{(p-\delta^2+\frac{2\chi}{\pi}\delta)^2+4\delta^2p}}{2}
	\right\}^\frac{1}{2}.
\end{equation}

We remark that in the large $\delta$ limit, i.e., $\delta\rightarrow\infty$, where the chemotactic response becomes negligible, the minimum propagation speed becomes $2\sqrt{p}$ irrespective of the value of $dp$ as is expected in the Fisher/KPP equation \cite{art:37F,art:37KPP,book:75AW}.
On the other hand, in the large stiffness limit $\delta\rightarrow 0$, the minimum propagation speed becomes $2\sqrt{p}-\chi$ irrespective of the value of $dp$, although the uniform saturated state $\rho=1$ is linearly unstable in this limit.
In the large stiffness limit, the effect of retraction due to the chemotactic response is maximized at the propagating front.

Even more interesting is that Eqs.~(\ref{eq_cl}) and (\ref{eq_minc}) imply a counter-intuitive phenomenon, i.e., the backward propagating wave with a negative propagation speed, where the population density initially saturated in the stable state may transit toward an unstable state in the local population dynamics.
This will be focused in this paper. 

\subsection{Unimodal traveling wave solution}\label{sec_analy}

As observed in the previous subsection, when the stiffness of chemotactic response is large, the propagation speed $c$ decreases because cells are attracted to migrate toward a higher-concentration region of chemoattractant in the propagating front.
Furthermore, the chemotaxis of cells may also create a peak in population density due to the chemotaxis.
In this subsection, we analyze analytically the existence of traveling waves with a single peak in the large-stiffness limit of chemotactic response $\delta^{-1}\rightarrow \infty$.
In this limit, we use the stiff flux function as
\begin{equation}\label{eq_sign}
	U_0[\partial_x \log S]=
	\left\{
		\begin{array}{cc}
			\chi &\mathrm{for} \quad \partial_x\log S > 0,\\
			-\chi &\mathrm{for}\quad \partial_x\log S<0.
		\end{array}
		\right .
	\end{equation}
We note that in the large-stiffness limit, the uniform saturated state $\rho=1$ is always linearly unstable so that stationary periodic patterns appear instead of the traveling wave in numerical computations.
Surprisingly, with the stiff flux Eq.~(\ref{eq_sign}), analytical unimodal traveling wave solution of Eqs.~(\ref{eq_KS})--(\ref{eq_S}), which are therefore certainly unstable, can be also computed explicitly.

\begin{figure}[tb]
	\includegraphics*[width=2.5in]{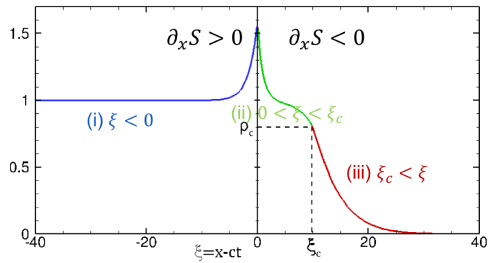}
	\caption{
		The schematic of domain decomposition in the relative coordinate system $\xi$ introduced in Sec.~\ref{sec_analy}.
	}\label{fig_domain}
\end{figure}

To obtain the analytical solution, we first suppose that $\tilde S(\xi)$ is a smooth unimodal function whose maximum is located at $\xi=0$, i.e.,
\begin{equation}\label{eq_dtS}
	\tilde S'(0)=0.
\end{equation}
Interestingly, with the piecewise linear property of the proliferation rate Eq.~(\ref{eq_P}) and the stiff flux function Eq.~(\ref{eq_sign}), the originally nonlinear FLKS system Eq.~(\ref{eq_KS}) is decomposed into the three linear equations in each different region as depicted in Fig.~\ref{fig_domain}, where the gradient of chemoattractant is positive $S'>0$ in region (i) ($\xi<0$) and negative $S'<0$ in regions (ii) and (iii) ($\xi> 0$). 
That is,
\begin{subequations}\label{eq_cauchy}
	\begin{equation}  
		(\chi-c)\tilde\rho_\mathrm{(i)}'(\xi)=\tilde\rho_\mathrm{(i)}''(\xi)+1-\tilde\rho_\mathrm{(i)}(\xi),
		\quad (\xi<0),
	\end{equation}
	\begin{equation} 
		-(\chi+c)\tilde\rho_\mathrm{(ii)}'(\xi)=\tilde\rho_\mathrm{(ii)}''(\xi)+1-\tilde\rho_\mathrm{(ii)}(\xi),
		\quad (0<\xi<\xi_c),
	\end{equation}
	\begin{equation} 
		-(\chi+c)\tilde\rho_\mathrm{(iii)}'(\xi)
		=\tilde \rho_\mathrm{(iii)}''(\xi)+p\tilde \rho_\mathrm{(iii)}(\xi),
		\quad (\xi\ge \xi_c),
	\end{equation}
\end{subequations}
where the subscripts (i), (ii), and (iii) are used to distinguish the different regions as depicted in Fig.~\ref{fig_domain} and $\xi_c$ is defined as $\tilde\rho(\xi_c)=\rho_c$.
Note that $\xi_c$ is a given parameter at this stage, but will be determined uniquely by Eq.~(\ref{eq_dtS}) later on.

The first-order derivative of $\tilde\rho(\xi)$ has a jump at $\xi=0$ because of the stiff flux function Eq.~(\ref{eq_sign}), i.e., 
\begin{equation}\label{eq_jump} 
	\tilde\rho'_\mathrm{(i)}(0)-\tilde\rho'_\mathrm{(ii)}(0)=2\chi\tilde\rho(0).
\end{equation}
Equation~(\ref{eq_cauchy}) can be solved analytically together with the boundary condition Eq.~(\ref{eq_bc_trho}), jump condition condition Eq.~(\ref{eq_jump}), and continuity and smoothness conditions at $\xi=0$ and $\xi=\xi_c$, i.e.,
\begin{equation}\label{eq_bc2_trho}
	\tilde \rho_\mathrm{(i)}(0)=\tilde\rho_\mathrm{(ii)}(0),
	\quad \tilde \rho_\mathrm{(ii)}(\xi_c)=\tilde\rho_\mathrm{(iii)}(\xi_c),
	\quad \tilde \rho_\mathrm{(ii)}'(\xi_c)=\tilde\rho_\mathrm{(iii)}'(\xi_c).
\end{equation}

From the dispersion relation between the decay rate at $\xi\gg 1$ and the traveling speed, Eq.~(\ref{eq_cl}), the traveling speed for the stiff flux function is written as
\begin{equation}\label{eq_sign_cl}
	c(\lambda)=\lambda+\frac{p}{\lambda}-\chi.
\end{equation}
Thus, the minimum speed $c_{\mathrm{min}}$ is obtained from the double root of the above equation as
\begin{equation}\label{eq_sign_minc}
	c_{\mathrm{min}}=2\sqrt{p}-\chi,
\end{equation}
with
\begin{equation}\label{eq_sign_lambda}
	\lambda=\sqrt{p}.
\end{equation}
In the following of this subsection, we seek for a unimodal traveling wave solution with the minimum traveling speed Eq.~(\ref{eq_sign_minc}).

By introducing a rescaled coordinate $\hat \xi$ defined as $\hat \xi=\sqrt{p}\xi$ and normalized parameters $s$ and $t$, $0<s,\,t<1$, defined as
\begin{equation}\label{eq_st}
	s=\exp(-\sqrt{1+p}\xi_c)=\exp(-\frac{\hat \xi_c}{t}),\quad t=\sqrt{1-\rho_c}=\sqrt{\frac{p}{1+p}},
\end{equation}
the solution of Eq.~(\ref{eq_cauchy}) together with the boundary conditions Eqs. (\ref{eq_bc_trho}), (\ref{eq_jump}), and (\ref{eq_bc2_trho}) is written explicitly as 
\begin{equation}\label{eq_trho1}
	\tilde\rho_{(i)}(\xi)=1+\alpha\exp(\nu\hat \xi),
\end{equation}
\begin{equation}\label{eq _trho2}  
	\tilde\rho_{(ii)}(\xi)=
	1-\beta\exp(\eta_+(\hat \xi-\hat \xi_c))
	+(\beta-t^2)\exp(\eta_-(\hat \xi-\hat \xi_c)),
\end{equation}
\begin{equation}\label{eq_trho3}  
	\tilde\rho_{(iii)}(\xi)=(1-t^2+\gamma(\hat \xi-\hat \xi_c))\exp(-(\hat \xi-\hat \xi_c)),
\end{equation}
where the decay rates $\nu$ and $\eta_\pm$ are written as
\begin{equation}\label{eq_nu}
	\nu=\hat\chi-1+\sqrt{\hat\chi^2-2\hat\chi+\frac{1}{t^2}}\quad (>0),
\end{equation}
\begin{equation}\label{eq_eta} 
	\eta_\pm=-1\pm\frac{1}{t},
\end{equation}
and the constants $\alpha$, $\beta$, and $\gamma$ are written as 
\begin{equation}\label{eq_alpha}
	\alpha=\frac{2\hat\chi(1-s^2)-2ts^{1-t}}
	{\mu_+s^2-\mu_-},
\end{equation} 
\begin{equation}\label{eq_beta}  
	\beta=\frac{2\hat \chi s^{1+t}-t^2\mu_-}{\mu_+s^2-\mu_-},
\end{equation}
\begin{equation}\label{eq_gamma} 
	\gamma=\frac{-4\hat\chi s^{1+t}+\mu_+s^2t(1+t-t^2)-\mu_-t(1-t-t^2)}
	{(\mu_+s^2-\mu_-)t}.
\end{equation}
Here, $\hat \chi$ is defined as 
\begin{equation}\label{tchi}
\hat\chi =\chi/\sqrt{p},
\end{equation}
and $\mu_\pm$ are defined as
\begin{align}\label{eq_mu}
	\mu_\pm&=\eta_\pm-\nu+2\hat\chi,\nonumber\\
	&=\hat\chi-\sqrt{\hat\chi^2-2\hat\chi+\frac{1}{t^2}}\pm\frac{1}{t}.
\end{align}

We note that $\mu_+$ is positive and $\mu_-$ is negative for any positive $\hat \chi>0$ because the derivative of $\mu_\pm$ with respect to $\hat \chi$ is always positive, i.e.,
\begin{equation} 
	\frac{d\mu_\pm}{d\hat\chi}
	=1-\frac{\hat\chi-1}{\sqrt{\hat \chi^2-2\hat\chi+\frac{1}{t^2}}}
	> 0,
\end{equation}
and $\mu_+$ is zero at $\hat \chi=0$ and $\mu_-$ is negative at $\hat\chi\rightarrow\infty$.

\begin{figure}[tb]
	\includegraphics*[width=4.5in]{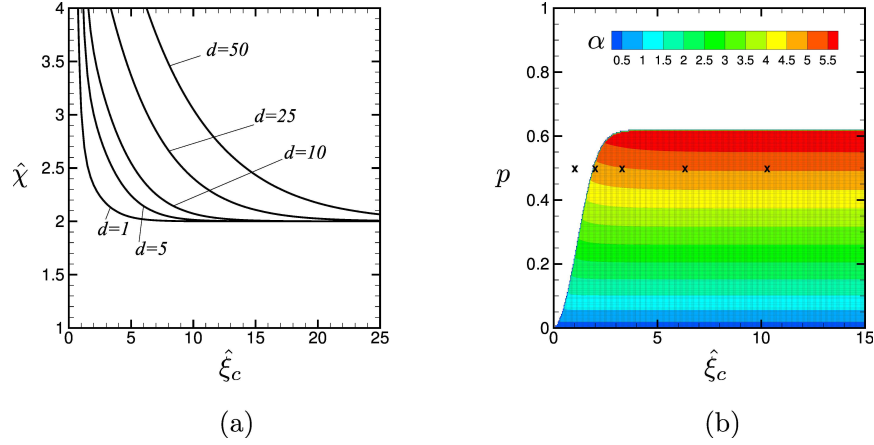}
	\caption{Figure (a) shows the solution curves of Eq.~(\ref{eq_F}) in $\hat \xi_c$--$\hat \chi$ plane with variation in the diffusion constant $d$, while the proliferation rate $p=0.5$ is fixed. 
	Figure (b) shows the parameter regime which satisfies the constraints Eqs.~(\ref{eq_const_f}) and (\ref{eq_const_g}).
	Here the modulation amplitude $\hat \chi=3.0$ is fixed.
	The contour shows the peak value of population wave, i.e., $\alpha$ defined in Eq.~(\ref{eq_alpha}).
	The symbols ``$\times$'' in figure (b) show the solutions of Eq.~(\ref{eq_F}) with $\hat \chi$=3.0 and $d=$1, 5, 10, 25, and 50, respectively, from left to right. 
}\label{fig_F_alpha}
\end{figure}

In order that Eqs. (\ref{eq_trho1})--(\ref{eq_mu}) constitute a unimodal traveling wave, the positivity of population-density gradient in $\xi<0$, i.e., $\alpha>0$, and the positivity of population density at $\xi\gg 1$, i.e., $\gamma>0$, must be satisfied.
These conditions give the following constraints among the parameters, i.e.,
\begin{equation}\label{eq_const_f}
	f(s,t)=\hat \chi(s^{t-1}-s^{t+1})-t >0,
\end{equation}
\begin{equation}\label{eq_const_g}
	g(s,t)=\mu_+s^2t(1+t-t^2)
	-\mu_-t(1-t-t^2)
	-4\hat\chi s^{1+t}>0.
\end{equation}
Figure \ref{fig_F_alpha}(b) shows a parameter regime of $\xi_c$ and $p$ which satisfies the constraints Eqs.~(\ref{eq_const_f}) and (\ref{eq_const_g}) simultaneously when the modulation amplitude $\hat\chi=3$ is fixed.
We note that Eqs.~(\ref{eq_const_f}) and (\ref{eq_const_g}) are independent on the diffusion coefficient $d$.
Remarkably, this shows the unimodal traveling wave solution exists in a certain range of parameters under the prerequisite assumption Eq.~(\ref{eq_dtS}).

We also remark that the upper bound of $p$ for large $\xi_c$, $p^u$ is calculated as $p^u=\frac{\sqrt{5}-1}{2}$ because Eq.~(\ref{eq_const_g}) for $\xi_c \gg 1$, i.e., $s\ll 1$, can be written as
\begin{equation}
	g(s,t)=-\mu_-t(1-t-t^2)-4\hat\chi s^{1+t}+\mathcal{O}(s^2).
\end{equation}
The second term is always negative and the first term is only positive for $t<t^u=\frac{\sqrt{5}-1}{2}$, which leads to $p<p^u$ as described above.

So far, we have constructed the unimodal traveling wave solution under the prerequisite assumption for the unimodal profile of chemoattractant $\tilde S(\xi)$ without considering the coupling of population density $\tilde \rho$ and chemoattractant $\tilde S$ via Eq.~(\ref{eq_S}).
We now consider the consistency of our analytical solution and the prerequisite assumption for $\tilde S(\xi)$ via Eq.~(\ref{eq_S}).

The unimodal profile of $\tilde \rho(\xi)$ warrants the unimodal profile of chemoattractant $\tilde S(\xi)$ via Eq.~(\ref{eq_S}), which is rigorously proved in Ref. \cite{art:16C}.
Thus, if and only if Eq.~(\ref{eq_dtS}) is satisfied our analytical solution Eqs.~(\ref{eq_trho1})--(\ref{eq_mu}) is proved to be an unimodal traveling wave solution of the FLKS system Eqs.~(\ref{eq_KS})--(\ref{eq_S}) for the stiff flux function Eq.~(\ref{eq_sign}).

From Eq.~(\ref{eq_tS}), the first derivative of $\tilde S(\xi)$ at $\xi=0$ is rewritten as
\begin{equation}\label{eq_dtS2}
	\tilde S'(0)=\frac{1}{2d}\int_0^\infty e^{-\frac{\zeta}{\sqrt{d}}}
	(\tilde \rho(\zeta)-\tilde \rho(-\zeta))d\zeta.
\end{equation}
By substituting Eqs. (\ref{eq_trho1})--(\ref{eq_trho3}) into Eq.~(\ref{eq_dtS2}), we can write Eq.~(\ref{eq_dtS}) as
\begin{equation}\label{eq_F} 
\begin{split}
	&F(\xi_c)=2d\sqrt{p}S'(0) \\
	&=-\alpha\left(\nu+\frac{1}{\sqrt{pd}}\right)^{-1}
	-e^{-\frac{\hat \xi_c}{\sqrt{pd}}} 
	\left\{	
		\left(1+\frac{1}{\sqrt{pd}}\right)^{-1}
		\left(
		\sqrt{pd}+t^2-\gamma\left(1+\frac{1}{\sqrt{pd}}\right)^{-1}
	\right) \right.\\
		&\left.	
		+\beta\frac{1-e^{-\left(\eta_+-\frac{1}{\sqrt{pd}}\right)\hat\xi_c}}{\eta_+-\frac{1}{\sqrt{pd}}}
	-(\beta-t^2)\frac{1-e^{-\left(
		\eta_--\frac{1}{\sqrt{pd}}\right)\hat\xi_c}}
	 {\eta_--\frac{1}{\sqrt{pd}}} 
\right\}=0.
\end{split}
\end{equation} 
Thus, the parameter $\hat \xi_c$ is determined by Eq.~(\ref{eq_F}) with a given parameter set of $\chi$, $p$, and $d$.
Furthermore, the obtained parameter $\hat \xi_c$ must also satisfy the constraints Eqs.~(\ref{eq_const_f}) and (\ref{eq_const_g}) because the normalized parameter $s$ defined by Eq.~(\ref{eq_st}) depends on $\hat \xi_c$.

Figure \ref{fig_F_alpha}(a) shows the solution curves of Eq.~(\ref{eq_F}) in $\hat \xi_c$--$\hat \chi$ plane with variation in the diffusion constant $d$, while the proliferation rate $p=0.5$ is fixed.
The parameter $\hat \xi_c$ grows as the diffusion coefficient $d$ increases.
The intersections of solution curves with different diffusion constants $d$ and a fixed modulation amplitude $\hat \chi=3.0$ are also shown in Fig. \ref{fig_F_alpha}(b).
It is seen that the parameter $\hat \xi_c$ determined by Eq.~(\ref{eq_F}) for $p=0.5$, $\hat\chi=3.0$, and $d\ge 5$ satisfies Eqs.~(\ref{eq_const_f}) and (\ref{eq_const_g}) simultaneously while that for $d=1$ does not satisfy the latter condition.  

We also remark that the solution curve of Eq.~(\ref{eq_F}) converges to $\hat \chi=2$ as $\hat \xi_c\rightarrow \infty$ because $F(\xi_c)$ is written, at $\xi_c=\infty$, as 
\begin{equation}\label{eq_F_asymp}
	F(\xi_c\rightarrow\infty)=
	\frac{2\hat\chi}{\mu_-}\left\{
		\left(\nu+\frac{1}{\sqrt{pd}}\right)^{-1}
		-\left(\frac{1}{\sqrt{pd}}-\eta_-\right)^{-1}
	\right\},
\end{equation}
where $\nu$ is the increasing function with respect to $\hat \chi$ and is equal to $-\eta_-$ at $\hat \chi=2$ (See also Eqs.~(\ref{eq_nu}) and (\ref{eq_eta})).

In conclusion of this subsection, the unimodal traveling wave solution written by Eq.~(\ref{eq_trho1})--(\ref{eq_mu}) exists in a certain range of parameters, where the parameters $\chi$, $p$, and $d$ satisfy Eqs. (\ref{eq_const_f}), (\ref{eq_const_g}), and (\ref{eq_F}) simultaneously.
In fact, from Fig.~\ref{fig_F_alpha}, we can find that Eqs.~(\ref{eq_const_f}), (\ref{eq_const_g}), and (\ref{eq_F}) are simultaneously satisfied when the modulation amplitude $\hat \chi$ and diffusion coefficient $d$ are sufficiently large, say $\hat \chi>2$ and $d>5$, for $p=0.5$.

\section{Numerical analysis}\label{sec_numeric}

\subsection{Numerical scheme}

We consider a one-dimensional interval $x=[0,L]$ and divide the interval into a uniform lattice-mesh system as $x_i=i\Delta x$ ($i=0,\cdots,I$), where $\Delta x$ is the mesh interval and $x_i$ represents the node of each mesh interval.
We solve Eq.~(\ref{eq_KS}) by using the following finite difference scheme,
\begin{equation}\label{eq_scheme}
	\begin{split}
	&\frac{\rho_i^{n+1}-\rho_i^n}{\Delta t}=\\
	&-\frac{1}{\Delta x}\left\{
		U_{i+1}^n\left(\frac{\rho_{i+1}^n+\rho_i^n}{2}\right)
		-U_{i}^n\left(\frac{\rho_{i}^n+\rho_{i-1}^n}{2}\right)
	\right\}
	+
	\frac{\rho_{i+1}^n-2\rho_i^n+\rho_{i-1}^n}{\Delta x^2}
	+P[\rho_i^n]\rho_i^{n+1},
\end{split}
\end{equation}
where the flux $U_i^n$ is calculated as
\begin{equation}\label{eq_Ui}
	U^n_i=U\left[\frac{\log S_i^n -\log S_{i-1}^n}{\Delta x}\right].
\end{equation}
Here $\rho_i^n$ is the average density in the mesh interval $x\in[x_i,x_{i+1}]$ at time $t=n\Delta t$, which is defined as $\rho_i^n=\frac{1}{\Delta x}\int_{x_i}^{x_{i+1}}\rho(n\Delta t,x)dx$, and $\Delta t$ is the time step size.
In Eq.~(\ref{eq_scheme}), $\rho_{-1}$ is replaced with $2-\rho_0$ at the left boundary, i.e., $x=0\,(i=0)$, and $\rho_I$ is replaced with $-\rho_{I-1}$ at the right boundary, i.e. $x=L\,(i=I)$.
Equation (\ref{eq_S}) is descretized on the same mesh intervals,
\begin{equation}\label{eq_schemeS}
	-d\frac{S_{i-1}^n-2S_i^n+S_{i+1}^n}{\Delta x^2}+S_i^n=\rho_i^n,
\end{equation}
and the same boundary conditions as $\rho_i$ are applied.
We solve Eq.~(\ref{eq_schemeS}) implicitly.

Numerical computations are performed for the scaled time and space variables defined as $\hat t=p t$ and $\hat x=\sqrt{p}x$, respectively.
The length of one-dimensional interval $\hat L\,(=\sqrt{p}L)\,=1000$, the number of mesh intervals $I=10000$ (i.e., the mesh interval $\hat \Delta x=0.1$), and the time step size $\hat \Delta t=\hat \Delta x^2/4$ are fixed.

The initial condition is set as 
\begin{equation}\label{eq_ini}
	\rho^0_i=S^0_i=\left\{
	\begin{array}{cc}
		1& (0\le\hat x_i\le \hat L_0) \\
		0& (\hat L_0\le \hat x_i \le \hat L)
	\end{array}
	\right.
\end{equation}
and $\hat L_0$ is set as $\hat L_0=100.0$ unless otherwise stated.

The accuracy of the numerical scheme is checked with respect to the propagation speed and decay rate of population density in the propagating front and the dispersion relation between them in Table \ref{t_error}.
The propagation speed $c^*$ is calibrated by tracing the position of a tip of propagating front $x^*$, which is defined as $\rho(x^*)=10^{-20}$, during time period $\hat t=[300.0,400.0]$.
We also measure the exponential decay $\lambda$ in Eq.~(\ref{eq_app_trho}) by applying $\lambda^*=\frac{\partial \log(\rho)}{\partial x}|_{x=x^*}$ to the numerical solutions.
The dispersion relation obtained by Eq.~(\ref{eq_cl}) with $\lambda^*$, $c(\lambda^*)$ is also compared to the measured propagation speed $c^*$.

\begin{table}
	\begin{tabular}{c c c c}
		\hline\hline
		$I_f$ -- $I_c$  
		& $|\frac{c^*_f-c^*_c}{c^*_f} |$ 
		& $|\frac{\lambda^*_f-\lambda^*_c}{\lambda^*_f}|$
		& $\frac{c^*_f-c(\lambda_f^*)}{c(\lambda_f^*)}$\\
		\hline
		10000 -- 5000 &$1.7\times 10^{-3}$& $2.8\times 10^{-3}$ &$1.1\times 10^{-3}$ \\
		20000 -- 10000 & $3.6\times 10^{-4}$& $5.6\times 10^{-4}$ &$8.1\times 10^{-4}$\\
		\hline\hline
	\end{tabular}
	\caption{The accuracy tests performed with different numbers of mesh interval $I$, i.e., $I$=5000, 10000, 20000. 
	The subscripts $f$ and $c$ represent the finer and coarser mesh systems, respectively.
The traveling speeds $c^*$ and exponential decay $\lambda^*$ are directly measured from the numerical solutions. 
}\label{t_error}
\end{table}

\subsection{Results}
Numerical simulations are performed for various values of modulation amplitude $\chi$ and stiffness $\delta$ while the diffusion coefficient $d=4.0$ and proliferation rate $p=0.5$ are fixed unless otherwise stated.
\begin{figure}[tb]
	\includegraphics*[width=2.5in]{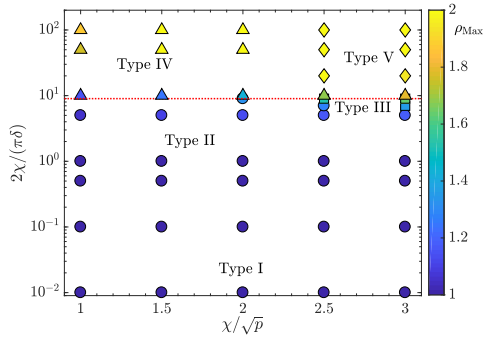}
	\caption{
		The diagram of different types of numerical solutions with variation in the modulation $\chi$ and stiffness $\delta^{-1}$.
		The circles (Type I and II) refer to Fig.~\ref{fig_forward}, the squares (Type III) refer to Fig.~\ref{fig_backward}, the triangles (Type IV) refer to Fig.~\ref{fig_periodic}, and the diamonds (Type V) refer to Fig.~\ref{fig_spik}.
		The colors of each symbol show the maximum value of population density in the spatial profile.
		The diffusion coefficient $d$ and proliferation rate $p$ are fixed as $d=4$ and $p=0.5$.
		The dotted horizontal line shows the critical value determined by the instability condition.
	}\label{fig_phase}
\end{figure}
\begin{figure}[tb]
	\includegraphics*[width=2.5in]{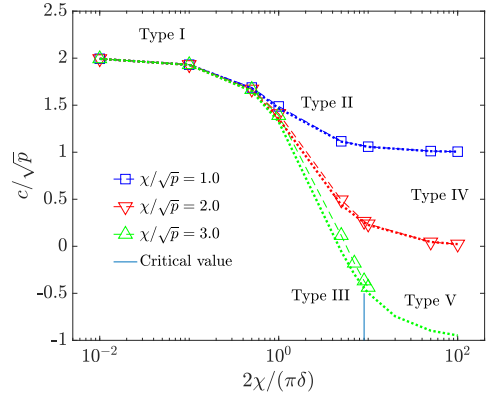}
	\caption{
		The speed of propagating front with variation in the modulation $\chi$ and stiffness $\delta^{-1}$.
		The values of parameters $d$ and $p$ are same as in Fig.~\ref{fig_phase}.
		The vertical line on the horizontal axis indicates the critical value of the instability condition.
		The dotted lines show the minimum speeds obtained by Eq.~(\ref{eq_minc}).
		The numbers (I)--(V) illustrate the solution types shown in Fig.~\ref{fig_phase}. 
	}\label{fig_speed}
\end{figure}
Figure \ref{fig_phase} shows the diagram of solution types obtained in the numerical simulations with variation in the modulation amplitude $\hat \chi$ and stiffness $2\chi/(\pi\delta)$.

\begin{figure}[tb]
	\includegraphics*[width=4.5in]{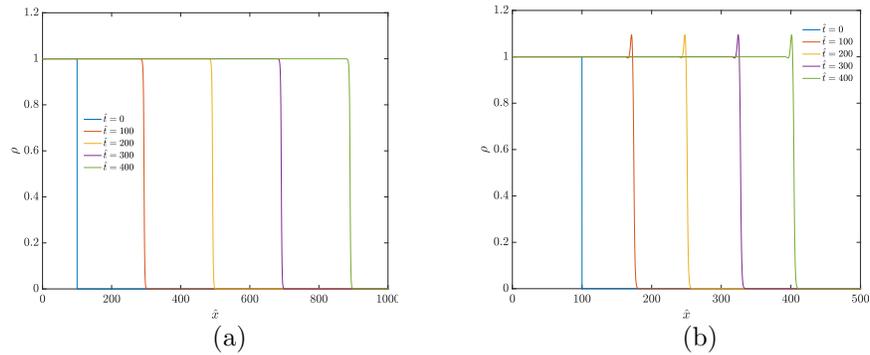}
	\caption{
		The snapshots of monotonic and non-monotonic traveling waves for $2\chi/(\pi\delta)=0.01$ (a) and $2\chi/(\pi\delta)=5.0$ (b), respectively. The other parameters are set as $\hat \chi=1.5$, $d=4.0$, and $p=0.5$. 
	}\label{fig_forward}
\end{figure}

\begin{figure}[tb]
	\includegraphics*[width=2.5in]{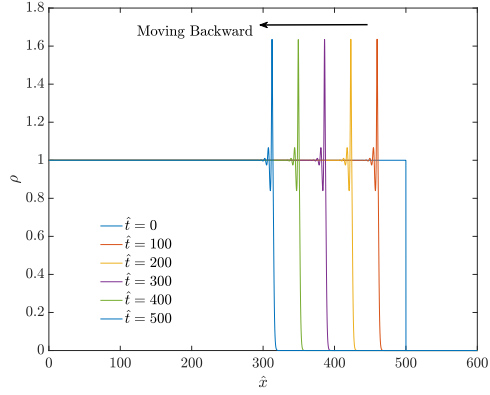}
	\caption{
		The snapshots of backward traveling wave for $\tilde \chi=3.0$, $2\chi/(\pi\delta)=9$, $d=4.0$, and $p=0.5$.
	}\label{fig_backward}
\end{figure}

\begin{figure}[tb]
	\includegraphics*[width=4.5in]{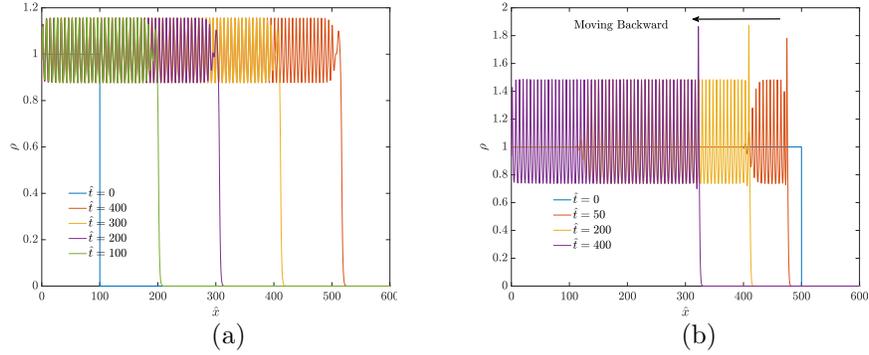}
	\caption{
		The snapshots of periodic pattern formation with a moving front in forward (a) and backward (b) directions. 
		The modulation parameter is set as $\hat\chi=1.0$ (a) and $\hat\chi=3.0$ (b). The other parameters are set as $2\chi/(\pi\delta)=10.0$, $d=4.0$, and $p=0.5$ in both (a) and (b).
	}\label{fig_periodic}
\end{figure}

Surprisingly, it is found that five different solution types exist, i.e., 
\begin{enumerate}[(I)]
	\item Monotonically decreasing traveling waves with a positive propagation speed. See Fig. \ref{fig_forward}(a).
	\item Non-monotonic traveling wave with a positive propagation speed. See Fig. \ref{fig_forward}(b).
	\item backward traveling wave, where the population wave with a steep peak propagates backward with a negative constant speed. See Fig. \ref{fig_backward}.
	\item Sequential strip pattern formation with a propagating front, where the stationary periodic pattern forms due to the instability around the uniform state $\rho=1$ while the front propagates with a constant positive or negative speed. See Fig. \ref{fig_periodic}
	\item Stationary localized spikes, which are synchronously created in the region of initially saturated state $\rho=1$.
\end{enumerate}
\begin{figure}[tb]
	\includegraphics*[width=2.5in]{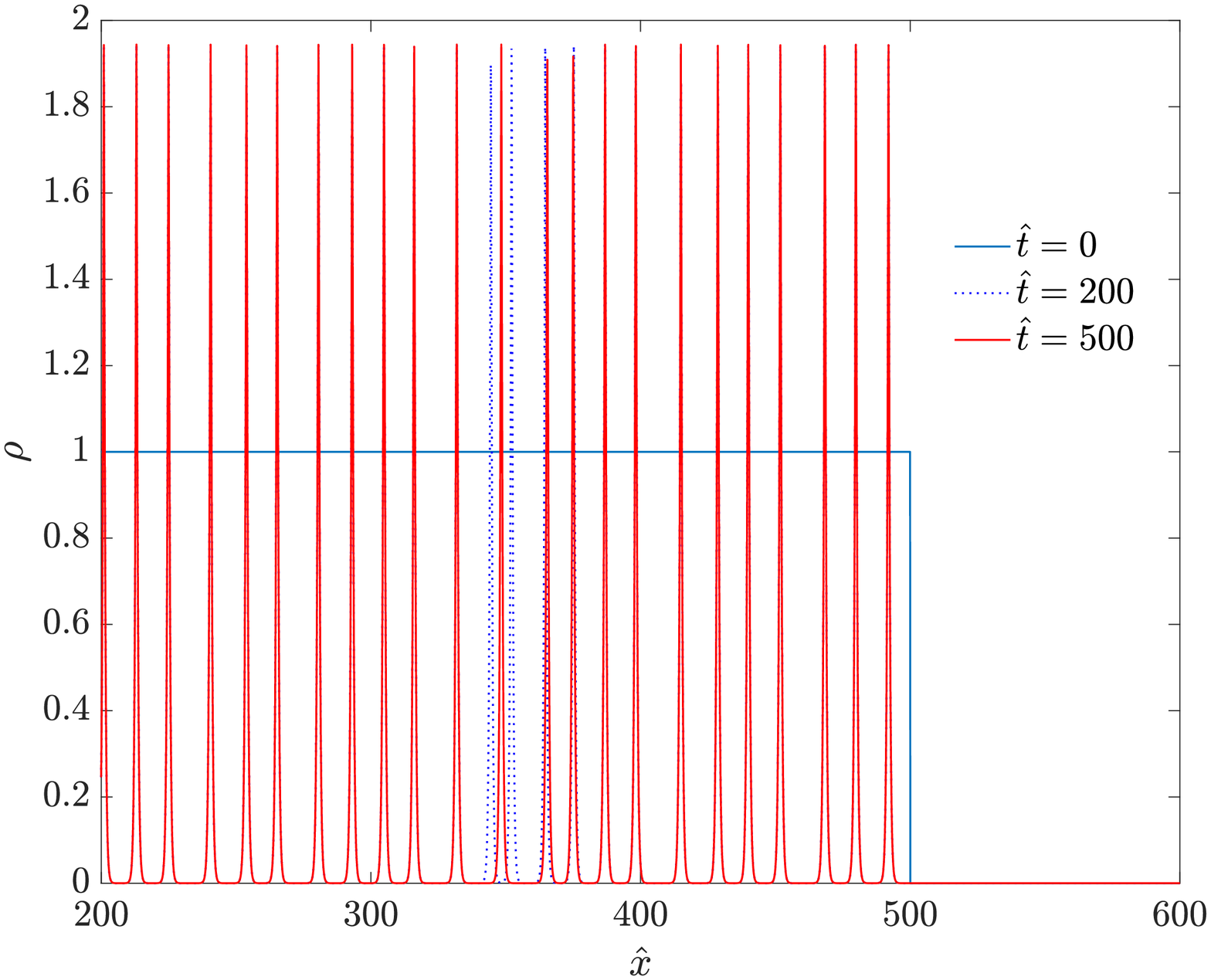}
	\caption{
		The snapshots of pattern formation of localized spikes. The parameters are set as $\hat \chi=3.0$, $2\chi/(\pi\delta)=20.0$, $d=4.0$, and $p=0.5$.
	}\label{fig_spik}
\end{figure}

Traveling waves, which propagate with constant speeds while keeping the spatial profiles, always appear under the linear stability condition Eq.~(\ref{eq_stability}).
The population density is concentrated due to the chemotaxis flux toward the gradient of chemoattractant produced by cells, so that non-monotonic traveling wave, which has a peak aggregation behind the propagating front, appears for a large stiffness parameter.
Surprisingly, backward traveling waves appear when modulation amplitude is sufficiently large, say $\hat\chi>2.0$, but the stiffness is slightly smaller than the critical value determined by the linear stability condition.
Interestingly, in the backward traveling waves, the local population initially saturated in the stable state $\rho=1$ transits toward the unstable state $\rho=0$ in the local population dynamics.
Furthermore, we can also see the transition of solution types from the backward traveling wave to the stationary localized spikes as increasing the stiffness parameter when the modulation amplitude is large $\hat\chi >2$.

Figure \ref{fig_speed} shows the propagation speed of the front which connects $\rho=0$ and $\rho=1$ with variation in the stiffness parameter $2\chi/(\pi\delta)$.
When the stiffness is sufficiently small, the propagation speed approaches to $\hat c=2.0$, which is the same as the traveling speed obtained in the Fisher/KPP equation without chemotaxis and coincides with the minimum speed Eq.~(\ref{eq_minc}) in the small stiffness limit.
As stiffness increases, the propagation speed decreases and the non-monotonic profile is created due to the chemotaxis.
Especially, when the modulation amplitude of chemotaxis $\chi$ surpasses the proliferation rate $p$ as $\chi/\sqrt{p}>2.0$, the retraction of the wave front occurs for a large stiffness.
The propagation speed converges to the minimum speed, Eq.~(\ref{eq_sign_minc}), in the large-stiffness limit when the modulation amplitude does not surpass the critical value mentioned above, i.e., $\hat\chi \le 2.0$.
However, when the modulation amplitude surpasses the critical value the propagation speed does not reach the minimum speed, instead the stationary localized spikes are created in the initially saturated state.

In Fig.~\ref{fig_speed}, we also plot the minimum traveling speed defined by Eq.~(\ref{eq_minc}).
Remarkably, except for the solution type V, the propagation speeds measured from the numerical results are close to the minimum speed Eq.~(\ref{eq_minc}).
The propagation speeds of numerical results seems to converge to the minimum speeds both in the small- and large-stiffness limits, i.e., $2\chi/(\pi\delta)\rightarrow$ 0 or $\infty$, while only small deviations are observed in intermediate regime.
The comparison of the propagation speeds of numerical results and the minimum traveling speed Eq.~(\ref{eq_minc}) is discussed in more detail in Sec. \ref{sec_compari_twspeed}.

\section{Discussion}
\subsection{The traveling speed}\label{sec_compari_twspeed}
\begin{figure}[tb] 
	\includegraphics*[width=2.5in]{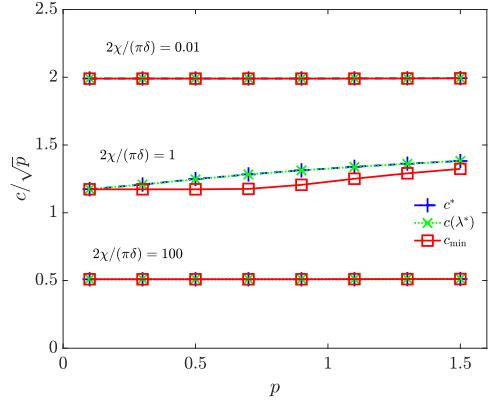}
	\caption{
		The comparison of the traveling speed measured from numerical solution $c^*$ to the minimum traveling speed $c_\mathrm{min}$ obtained by Eq.~(\ref{eq_minc}).
		The modulation amplitude $\hat\chi=1.5$ and diffusion constant $d=1.0$ are fixed while the proliferation rate $p$ varies.
		The dispersion relation Eq.~(\ref{eq_cl}) between the traveling speed and exponential decay of population density far ahead the front, i.e., $c(\lambda^*)$ is also plotted. 
	}\label{fig_cplot}
\end{figure}

\begin{figure}[tb]
	\includegraphics*[width=2.5in]{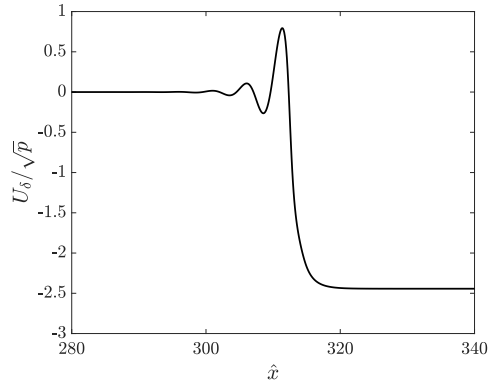}
	\caption{
		The snapshot of chemotactic drift speed $U_\delta$ for Fig. \ref{fig_backward} (i.e., Type III in Fig.~\ref{fig_speed}) at time $\hat t=500$. 
	}\label{fig_U}
\end{figure}

We observed a wide range of parameters for which traveling waves are propagating in the numerical scheme.
We found numerical evidence that the dispersion relation Eq.~\eqref{eq_cl} between the propagation speed $c(\lambda)$ and the exponential decay rate $\lambda$ is satisfied far ahead the front ($x\gg ct$). See Fig.~\ref{fig_cplot}.
However, this does not allow to compute readily the actual speed of propagation $c$, as it is the case for the so-called {\em pulled fronts} in reaction-diffusion equations. 
The notion of pulled front corresponds to those reaction-diffusion traveling waves for which the dynamics of small densities drive the whole expansion of the population. 
In particular, the remarkable formula $c = \min_\lambda c(\lambda)$ holds in such regimes. 
A celebrated example is the Fisher/KPP equation, or more generally any classical reaction-diffusion for which the maximal growth rate per capita is reached at zero density of individuals, {\em e.g.}
\begin{equation}
\partial_t \rho = \partial_{xx} \rho + P[\rho]\rho\,, \quad \max_{\rho\geq 0} P[\rho] = P[0]\, .
\end{equation}
This is usually opposed to the notion of {\em pushed fronts}, for which the whole range of individuals contribute to the expansion dynamics. 
As such, there is generally no explicit formula available for the speed.

Here, the discrepancy between $c$ and $\min_\lambda c(\lambda)$, even relatively small (see Figs. \ref{fig_speed} and \ref{fig_cplot}), is the signature of non-local effects, which shape the dynamics of expansion as in {\em pushed fronts}. 
Before discussing these non-local effects, it is noticeable that $c\sim\min_\lambda c(\lambda)$ actually holds true both as $\delta\to +\infty$ (no chemotaxis), and $\delta \to 0$ (stiff chemotactic response). 
The former is nothing but the reaction-diffusion limit, where the Fisher/KPP equation is recovered. 
In the latter case, the chemotactic drift converges towards a stepwise function taking value $\pm \chi$. 
In particular, it is expected that the drift has constant value $-\chi$ for $x$ far ahead. 
The situation is equivalent to a shifted Fisher/KPP equation on a (possibly moving) half-line. 
Hence, the formula for the speed $2\sqrt{p}-\chi$.

Now, we discuss the intermediate situation based on numerical insights. 
The intuition for pulled fronts is that propagation is a combination of diffusion and growth. 
The latter achieves its maximal rate in the tail of the population distribution far ahead, where the equation can be approximated with a linear equation which possesses explicit solutions travelling at the minimal speed. 
If the chemotactic transport speed $U_\delta$ would be maximal as $x\to +\infty$ (meaning minimal chemotactic effect, recall that $U_\delta<0$ far ahead), then we would expect the same conclusion as for pulled fronts. 
However, we observe the opposite (see Fig.~\ref{fig_U}), $U_\delta$ is indeed minimal as $x\to +\infty$ (meaning maximal chemotactic effect).

This maximal retreating chemotactic effect can balance the leading driving effect at small densities. 
This yields a wave traveling faster than the minimal speed, a signature of pushed fronts.

This conclusion is questionable, as it might be unrealistic to observe maximal chemotactic effect at lower density. 
However, we argue that the logarithmic sensing assumed in this model is indeed a way for bacteria to navigate across several order of magnitudes of chemical concentrations, and to modulate their response accordingly. 
Within this perspective, we observe this counter-balancing effect of chemotaxis drift vs. reaction-diffusion through the mechanism of expansion ({\em pushed} vs. {\em pulled}). 
This formal reasoning deserves more mathematical analysis beyond the numerical evidence presented in this work.

\subsection{Comparison to the unimodal analytical solution}
In Sec. \ref{sec_analy}, we show that the unimodal traveling wave solution of Eqs.~(\ref{eq_KS})--(\ref{eq_S}) with the stiff flux Eq.~(\ref{eq_sign}), i.e., $\delta\rightarrow 0$, which is certainly unstable, can be analytically computed in a certain parameter regime, e.g., $\hat\chi>2$, $p<\frac{\sqrt{5}-1}{2}$, and a sufficiently large $d$.
In this subsection, we compare the analytical solution with the numerical results for large stiffness parameters and discuss how the analytical solutions coincide or differ from the numerical solutions.

\begin{figure}[tb]
	\includegraphics*[width=4.5in]{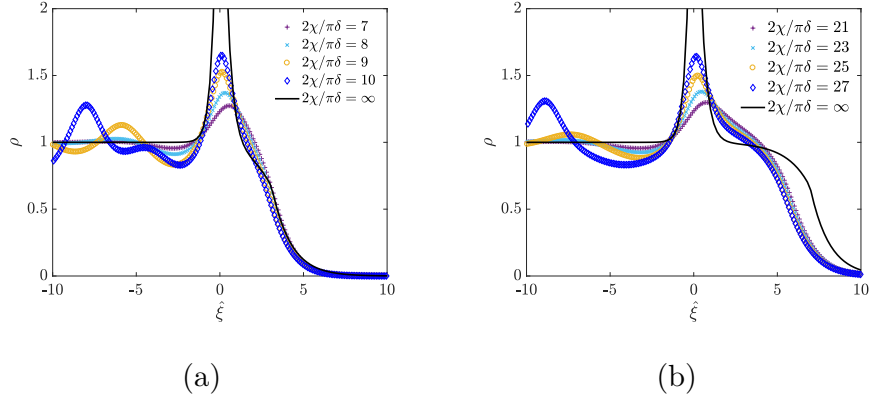}
	\caption{
		Numerical solutions for large-stiffness parameters are compared to the analytical solution for the stiff flux Eq.~(\ref{eq_sign}), i.e., $2\chi/(\pi\delta)\rightarrow \infty$.
		The modulation amplitude $\chi$ and proliferation rate $p$ are fixed as $\hat\chi=2.5$ and $p=0.5$, respectively. The diffusion coefficient $d$ is set as $d=4$ in figure (a) and $d=16$ in figure (b).
	}\label{fig_compari_analy}
\end{figure}

In Fig. \ref{fig_compari_analy}, population density profiles obtained for large stiffness are compared with the analytical solution in the stiff flux.
In Table \ref{table_xc}, decay rates $\lambda$ at $\xi\gg 1$ and the distances to the position for $\rho(\xi_c)=\rho_c$, $\xi_c$ are calculated from the profiles in Fig. \ref{fig_compari_analy}.

As increasing the stiffness parameter, the peak profile around $\hat \xi=0$ of the numerical solution is sharpened and the peak position approaches to $\hat \xi=0$.
The decay rate at $\hat \xi\gg 1$ also approaches to that of the analytical solution as increasing the stiffness parameter.

However, the nonmonotonic profiles of numerical solutions always oscillates behind the peak of propagation front, while the analytical solution of stiff flux is monotonic for $\hat \xi<0$.
Furthermore, the oscillation mode of numerical solution grows as increasing the stiffness parameter, and it bifurcates to the stationary oscillation from the traveling wave when the stiffness becomes larger than the critical value of the linear stability condition.
Thus, the unimodal traveling wave obtained analytically for the stiff flux does not appear, instead the nonmonotonic traveling wave which connects $\rho=0$ at $\xi \gg 1$ and oscillation in $\xi<0$ appears for a large stiffness parameter under the linear stability condition.

\begin{table} 
	\begin{tabular}{ccc c ccc}
		\hline\hline
		\multicolumn{3}{c}{$d=4.0$}&&\multicolumn{3}{c}{$d=16.0$}\\
		\cline{1-3}\cline{5-7}
		$\frac{2\chi}{\pi\delta}$&$\lambda/\sqrt{p}$&$\hat \xi_c$
		&&$\frac{2\chi}{\pi\delta}$&$\lambda/\sqrt{p}$&$\hat \xi_c$\\
		\hline
		7.0&1.30   & 3.08  && 21.0& 1.065 & 5.57\\
		8.0&1.26   & 2.96  && 23.0& 1.062 & 5.44\\
		9.0&1.23   & 2.79  && 25.0& 1.059 & 5.28\\
		10.0&1.20  & 2.61  && 20.0& 1.058 & 5.16\\
		$\infty$&1.00 & 3.09  && $\infty$&1.00 &6.95\\
		\hline\hline
	\end{tabular}
	\caption{The decay rate $\lambda$ defined in Eq.~(\ref{eq_app_trho}) and the distance from the peak of chemoattractant to the position where the population density equals to $\rho_c$, i.e., $\xi_c=x_c-x_S$ where $\rho(x_c)=\rho_c$ and $\partial_x S(x_S)=0$, with variation in the stiffness. The modulation amplitude $\chi$ and proliferation rate $p$ are fixed as $\hat \chi=2.5$ and $p=0.5$, respectively.}\label{table_xc}
\end{table}

\section{Concluding remarks}

Traveling wave and aggregation in the flux-limited Keller-Segel system Eqs.~(\ref{eq_KS})--(\ref{eq_S}), which describes the stiff and bounded chemotaxis flux to the logarithmic sensing of chemical cues, are investigated theoretically and numerically.

The numerical simulations uncover the variety of solution types in the FLKS system, i.e., (i) monotonic traveling wave for a small stiffness, (ii) nonmonotonic traveling wave with a positive propagation speed for a small modulation amplitude, i.e., $\hat \chi \le 2$, and a sufficiently large stiffness under the linear stability condition, (iii) backward traveling wave for a large modulation amplitude, i.e., $\hat\chi> 2$, and a sufficiently large stiffness under the linear stability condition, (iv) sequential strip pattern formation with a propagating front for a small modulation amplitude, i.e., $\hat\chi \le 2$, in the linear unstable condition, and (v) the stationary localized spikes for a large modulation amplitude, i.e., $\hat\chi > 2$, in the linear unstable condition.

Our study leads to several counter-intuitive outcomes.
Firstly, the non-monotonic backward traveling wave, with a complex profile, appear for a certain range of parameters, where a local population density initially saturated in the stable state $\rho=1$ transits towards the unstable state $\rho=0$ in the local population dynamics.
Secondly, transition occurs from retreating wave to stationary localized spikes, which are synchronously created only in the region of initially saturated state $\rho=1$, as increasing stiffness.
These behaviors stem from a large chemotaxis flux when the modulation amplitude is large as $\hat \chi>2$.

In the theoretical part, we obtain a novel analytic formula of the minimum propagation speed, Eq.~(\ref{eq_minc}), for a specific flux function Eq.~(\ref{eq_arcU}) from the general dispersion relation between propagation speed and exponential decay rate in the propagation front, Eq.~(\ref{eq_cl}).
Remarkably, except for localized spiky solutions, the traveling speeds of numerical results are asymptotically close to the minimum propagation speed both in the small- and large-stiffness limits, while they slightly deviate from the minimum propagation speed in the intermediate stiffness regime due to the counter-balancing effect of chemotactic drift vs. reaction-diffusion through the expansion mechanism ({\em pushed} vs. {\em pulled}).

We also discover an analytical solution of unimodal traveling wave of the FKLS system with the stiff flux Eq.~(\ref{eq_sign}), although the solution is certainly unstable in this limit and thus does not appear in the numerical simulations.

Because of complexity, biological processes of pattern formation under the effect of chemotaxis in cells have yet to be elucidated.
However, interestingly, the present study demonstrates that the FLKS model can reproduce the sequential strip pattern formation in expanding population which is observed in bacterial experiments \cite{art:66A, art:11Letal}.
The localized spikes pattern obtained in our simulation may also resemble the spot array formation of chemotactic bacteria, which appears due to an active accumulation \cite{art:91BB}.
These results advocate the usage of the FLKS model in an elucidation of some aspects of the pattern formation mechanism of chemotactic cells.

\section*{Acknowledgements}

VC has received funding from the European Research Council (ERC) under the European Union's Horizon 2020 research and innovation program (grant agreement No 639638).
BP has received funding from the European Research Council (ERC) under the European Union's Horizon 2020 research and innovation programmed (grant agreement NO 740623).
SY has received funding from the Japan Society for the Promotion of Science (JSPS) KAKENHI Grant Numbers 15KT0110 and 16K17554.
Part of this work was done during the trimester ``Stochastic Dynamics Out of Equilibrium'' at the Institut Henri Porincar\'e -- Centre Emille Borel.
The authors also thank this institution for hospitality and support.

\medskip
Received xxxx 20xx; revised xxxx 20xx.
\medskip

\end{document}